\documentclass[11pt]{article}

\usepackage{epsf}
\usepackage{psfrag}
\usepackage{amsmath}
\usepackage{amssymb}
\usepackage{lscape}
\usepackage{latexsym}
\usepackage[usenames]{color}
\usepackage[mathcal]{eucal}
\usepackage{stmaryrd}
\usepackage{textcomp}
\usepackage{setspace} 
\usepackage{graphicx} 
\usepackage{pgfplots} 
\usepackage{epsfig}   
\usepackage{epstopdf} 
\usepackage{appendix}
\usepackage{afterpage}
\usepackage{rotating} 
\usepackage{float} 
\usepackage[misc,geometry]{ifsym} 
\usepackage{comment}

\newtheorem{theorem}{Theorem}[section]

\newtheorem{proposition}{Proposition}[section]

\newtheorem{example}{Example}[section]

\newcommand{\ds}{\displaystyle}
\newcommand{\dZ}{{\cal Z \kern -0.7em Z}}
\newcommand{\dC}{{\rm\hbox{C \kern-0.8em\raise0.2ex\hbox{\vrule height5.4pt width0.7pt}}}}
\newcommand{\dQ}{{\rm\hbox{Q \kern-0.85em\raise0.25ex\hbox{\vrule height5.4pt width0.7pt}}}}

\newcommand{\proofbox}{\hspace{\fill}{$\Box$}}

\newenvironment{proof}{Proof.}{\proofbox}

\newcommand\old[1]{}
\newcommand{\beqa}{\begin{eqnarray*}}
\newcommand{\eeqa}{\end{eqnarray*}}
\newcommand{\rectangle}{\fboxsep0pt\fbox{\rule{1em}{0pt}\rule{0pt}{1ex}}}

\begin{document}

\title{\bf \Large{Green Inventory Management: Leveraging Multiobjective Reverse Logistics}}

 \author{I. B. Wadhawan \thanks{School of Mathematical Science, Queensland University of Technology (QUT), 
 {\tt email: indu.wadhawan@qut.edu.au (corresponding author).}}
 \and
 M. M. Rizvi \thanks{Centre for Smart Analytics (CSA), Institute of Innovation, Science and Sustainability, Federation University Australia, {\tt email:  rizmm001@mymail.unisa.edu.au}}
  }
\maketitle
\textbf{Abstract} The paper proposes a novel Economic Production Quantity (EPQ) inventory model within a reverse logistics framework, addressing new and repaired products with varying quality and demand patterns. The model integrates production and remanufacturing rates as functions of lot sizes and cycle numbers to develop a feasible inventory cost function. A key contribution of the study is formulating a multiobjective optimization framework that simultaneously minimizes inventory costs and accounts for environmental sustainability by considering greenhouse gas (GHG) emissions and energy consumption during production processes. The problem is formulated as a mixed-integer nonlinear programming (MINLP) model, with integer constraints on lot sizes and cycle counts and a continuous return rate. Numerical case studies taking test problems from existing literature are used to validate the model through extensive sensitivity analyses. Both mathematical optimization and heuristic optimization methods are applied to solve multiobjective optimization problems, and Pareto fronts are illustrated along with the interpretation of the results. The results, obtained using solvers in MATLAB and AMPL, highlight the model’s ability to balance operational efficiency and environmental responsibility. Pareto frontiers generated from the analysis provide strategic insights for decision-makers seeking to optimize cost and sustainability in inventory systems.


\textbf{AMS} subject classifications. 90B05, 90C25, 90C29, 90C30.


\section{Introduction}\label{sec1}
\noindent Reverse logistics (RL) plays a crucial role in modern supply chain management by enabling the movement of products from consumers back to manufacturers for return, repair, refurbishment, recycling, or proper disposal. This process complements the traditional forward flow, where goods travel from manufacturers to end users. While forward logistics focuses on product delivery and customer satisfaction, reverse logistics adds value by recovering product worth, reducing environmental impact, and minimizing resource consumption \cite{Schrady1967,Dobos2004, Dobos2006}. The growing importance of RL lies in its ability to extend product life cycles, reduce operational costs, and support sustainable business practices. As environmental regulations are tightened and consumer awareness increases, reverse logistics has become an essential strategy for companies aiming to improve efficiency, reduce waste, and maintain a competitive edge through responsible resource management.

\noindent In early work, Schrady \cite{Schrady1967} presented a deterministic Economic Order Quantity (EOQ) model for repaired items, and Nahmias and Rivera \cite{Nahmias1979} later extended it by allowing a finite repair rate within a production cycle. 
Richter \cite{Richter1996a, Richter1996b} introduced an Economic Order Quantity (EOQ) framework to determine the optimal batch sizes for production and remanufacturing, considering the product return rate. In a later study \cite{Richter1997}, the author investigates a system with two inventories and managed under either a “dispose-all” or “recover-all” policy.
 Dobos and Richter \cite{Dobos2004} extended this by considering limited production and repair rates and later included the possibility that some returned products might not be suitable for recycling \cite{Dobos2006}.  El~Saadany and Jaber \cite{Saadany2010} suggested a mixed strategy of production and remanufacturing, while El~Saadany \cite{Saadany2013} further explored reverse logistics by distinguishing unlimited remanufacturing from normal production and adding disposal costs. Environmental aspects are introduced by Jaber et~al.\ \cite{Jaber2013a} with a two-level vendor–buyer model that includes greenhouse gas (GHG) emissions. Singh and Sharma\ \cite{Singh2016} later incorporated finite production and remanufacturing rates considering EPQ model. Bazan et~al.\ \cite{Bazan2015, Bazan2016} proposed a four-objective model addressing waste, GHG emissions, and energy use. Kozlovskaya et~al.\ \cite{Kozlovskaya2017} studied deterministic inventory systems with constant demand and returns to minimize production, holding, and switching costs, and Singh et~al.\ \cite{Singh2020} considered variable production and remanufacturing under learning effects.

\noindent Earlier studies have shown that researchers consider remanufactured products to be as good as new ones. In real life, many customers view repaired or recycled items as inferior quality. Jaber and El Saadany \cite{Jaber2009}, Hasanov et al. \cite{Hasanov2012}, and Jaber and Peltokorpi \cite{Jaber2024} present models that create two separate markets for new and repaired products for EOQ models. Later, Singh and Sharma \cite{Singh2016} and Sharma et al. \cite{Swati2021} introduced models where production and remanufacturing rates are finite, and they developed EPQ models for reverse logistics systems. However, these models consider production rates and remanufacturing rates as fixed constants, with lot sizes as continuous variables.

\color{black}

\noindent This study addresses the gap in the literature in three main ways:
\begin{itemize}
    \item Integrated production and remanufacturing rates as variables: our proposed models integrate production and remanufacturing rates as functions of lot sizes and number of cycles for new and remanufactured items.
    \item Integer decisions variables: we consider lot sizes and cycle counts to be integers, as they are in real production, turning the problem into a mixed-integer programming model.
    \item Environmental goals: we add a multiobjective framework that  minimises not only inventory cost but also greenhouse gas emissions and energy use.
\end{itemize}

\noindent Therefore, the research goals of this analysis are:
\begin{enumerate}
    \item Create a material flow model for reverse logistics that finds the lowest possible inventory cost.
    \item Develop mathematical expressions to identify the best production rates, lot sizes, and number of cycles that keep inventory costs minimum.
    \item Reduce environmental impact by including greenhouse gas (GHG) emissions and energy use in both production and remanufacturing within the reverse logistics system.
    \item Develop multiobjective models that measure both inventory cost and environmental cost, and minimize these competing goals at the same time using different mathematical techniques.
    \item Design and test algorithms to solve the proposed models and verify their performance through computer experiments.
\end{enumerate}

\noindent Unlike existing reverse logistics EPQ models that assume fixed production and remanufacturing rates, the proposed formulation endogenizes both rates as functions of lot sizes and cycle numbers. This modelling feature creates a direct interaction between operational decisions and environmental performance measures. Furthermore, the incorporation of integer lot sizes and cycle frequencies transforms the problem into a mixed-integer nonlinear programming framework that more closely reflects practical production planning environments. While the model builds upon established reverse logistics structures, it extends the literature by integrating operational and environmental decision variables within a unified optimization framework capable of generating implementable inventory policies.

\noindent The rest of this paper is organised as follows: Section \ref{Insight} explains the main theoretical and practical contributions of the research. Section \ref{Form} describes all the variables used in this paper. Section \ref{mathmodel} introduces a new mathematical model for an EPQ reverse logistics system. It also presents numerical experiments and sensitivity analyses for mixed-integer optimization problems aimed at reducing inventory costs. Section \ref{EnvImp} describes the greenhouse gas (GHG) emissions and energy use linked to both production and remanufacturing. Section \ref{MOP} develops a two-objective optimization model that looks at both key factors: inventory costs and GHG emissions from production and remanufacturing. Detailed formulas for the cost function  are provided in the appendix.

\section{Decision and Managerial Insights} \label{Insight}
\noindent The proposed models comprise two main functional components. The first focuses on inventory management, to minimise inventory costs as described in Section \ref{mathmodel}. In this context, supply chain managers and decision-makers can use Model \eqref{eqmathemodel1} to determine production and remanufacturing rates, integer lot sizes, cycle numbers, and the quality level of returned items from repaired products ($0 \leq s \leq 1$) originating from the secondary market.
The second objective, discussed in Section \ref{MOP}, addresses multiobjective optimization for environmental sustainability, aiming the reduction of greenhouse gas (GHG) emissions and energy consumption. Decision-makers in manufacturing, remanufacturing, sustainability initiatives, and circular economy practices can apply Models \eqref{modelmop1} and \eqref{modelmop2} to evaluate trade-offs between inventory costs and environmental impacts, enabling strategies that reduce GHG emissions and energy use while protecting the environment.
The inventory management framework can also be adapted for manufacturing companies by modifying constraints, for example, maintaining production control by fixing lot sizes while still meeting customer demand.

\noindent The proposed models incorporate varying demand levels for both new and repaired products while addressing environmental pollution-related costs, a factor often neglected in conventional inventory management models. We develop multiobjective mathematical formulations to examine the trade-offs between inventory costs and costs associated with GHG emissions, as well as between inventory costs and energy consumption costs. By integrating a multiobjective approach with scalarization techniques, the model enables a systematic exploration of these trade-offs, balancing economic efficiency with environmental sustainability. This provides decision-makers with a more comprehensive, flexible, and environmentally conscious framework than traditional single-objective inventory models.

\noindent Implementing our model delivers direct environmental benefits by reducing greenhouse gas (GHG) emissions and lowering energy consumption through optimized production and remanufacturing processes. By explicitly integrating these environmental factors into inventory decision-making, the model helps companies minimize their carbon footprint and resource use, contributing to cleaner manufacturing and reverse logistics operations. 

\noindent The model is relatively straightforward to implement for companies where the demands are known in advance. This can be done by establishing data collection processes for inventory, demand forecasting, and cost analysis. Businesses that already use optimization tools or software platforms for logistics can easily integrate our mixed integer programming models into their existing systems. Moreover, the algorithm for finding Pareto solutions is designed to be adaptable to various real-world constraints.

\noindent We look at cases where new and refurbished products are viewed differently in quality. Take the tyre market as an example. New tyres and retreaded or repaired tyres serve different kinds of customers, so it makes sense to treat them as two separate markets. The first market sells new tyres to people who want the new product, while the second market offers repaired tyres at a lower price for cost-conscious buyers. Our model is based on this idea, as shown in Figure \ref{figflowModel1}. Making and remanufacturing products also affect greenhouse-gas emissions and energy use, so our approach aims to balance production costs with environmental impact. A similar case appears in the electronics sector, where refurbished phones, laptops, and tablets have steady demand. Another good example is the furniture industry, where restored wooden tables, chairs, and cabinets are sold alongside new pieces, requiring careful planning and inventory control that reinforcement learning can support.

\section{Formulations and Solutions of the Proposed Models }\label{Form}

\noindent Decision variables:
\begin{table} [H]
\begin{tabular}{l l}
$m$& number of remanufacturing cycles.\\
$n$& number of manufacturing cycles.\\
$Q_p$& production batch size (units).\\
$Q_r$& remanufacturing batch size (units).\\
$s$  & quality level of returns items from repaired items ($0\leq s \leq 1$). \\

\end{tabular}
\end{table}\newpage
\noindent Input parameters:
\begin{table} [H]
\begin{tabular}{l l}
$D_p$& demand rate for new items (units/unit-time)\\
$D_r$& demand rate for repaired items (units/unit-time)\\
$p$  & collection proportion of new items $0<p\leq 1$\\
$q$  & quality level of returns items from new items ($0\leq q \leq 1$) \\
$r$  & collection proportion of repaired items $0\leq r\leq 1$\\
$qpD_p+srD_r$& repairable stock rate ($R_1 = qpD_p, R_2 = srD_r$) (units/unit-time)\\
$S_p$ & setup cost of supply depot ($\$$/setup)\\
$S_r$ & setup cost of repair depot ($\$$/setup)\\
$h_p$ & holding cost of supply depot ($\$$/unit/unit-time)\\
$h_r$ & holding cost of repair depot ($\$$/unit/unit-time).\\
$a_p$, $a_r$ & emissions function parameter for new and repair items (ton year\(^2\)/unit\(^3\))\\
$b_p$, $b_r$ & emissions function parameter for new and repair items (ton year/unit\(^2\))\\
$c_p$, $c_r$ & emissions function parameter for new and repair items (ton year/unit)\\
$K_p$, $K_r$ & energy required by the production and repair machine to produce \\ & one unit (kWh/unit)\\
$W_p$, $W_r$ & idle power of the production and repair machine (kW).\\
$C_{GHG}$, $C_{ENG}$ & Cost per unit associated with GHG emission and energy consumption.
\end{tabular}
\end{table}

\noindent Decision variables dependent parameters:
\begin{table} [H]
\begin{tabular}{l l}
$T_p$& length of new items cycle, where $T_p=nQ_p/D_p$\\
$T_r$& length of repaired items cycle, where $T_r=mQ_r/D_r$\\
$T$& length of new procurement and repair cycle, where $T=T_p+T_r$\\
$P$ & production rate, where \(\ds P(Q_p, n) = \frac{D_p}{1 - \frac{2S_pD_p}{h_p(nQ_p)^2}}\), \(1 - \frac{2S_pD_p}{h_p(nQ_p)^2}>0 \)\\
$\lambda$ & remanufacturing rate, where \(\ds \lambda(Q_r, m) = \frac{D_r}{1 - \frac{2S_rD_r}{h_r(mQ_r)^2}}\), \(1 - \frac{2S_rD_r}{h_r(mQ_r)^2}>0\).
\end{tabular}
\end{table}
\par
\noindent In developing the proposed model, several key assumptions are made to reflect the practical complexities of managing new and repaired inventory. First, it is assumed that repaired items are of a different quality than new items, leading to differing demand patterns (see, \cite{Singh2016, Jaber2009}). We also assume that in the new procurement period $T_p$, the demand fulfils $D_p \geq nQ_p$. Similarly, in the repair cycle $T_r$,  $D_r \geq mQ_r$. 

\section{Mathematical Models and Solutions}\label{mathmodel}
\begin{figure}[H]
	\centering
\includegraphics[width=0.95\textwidth]{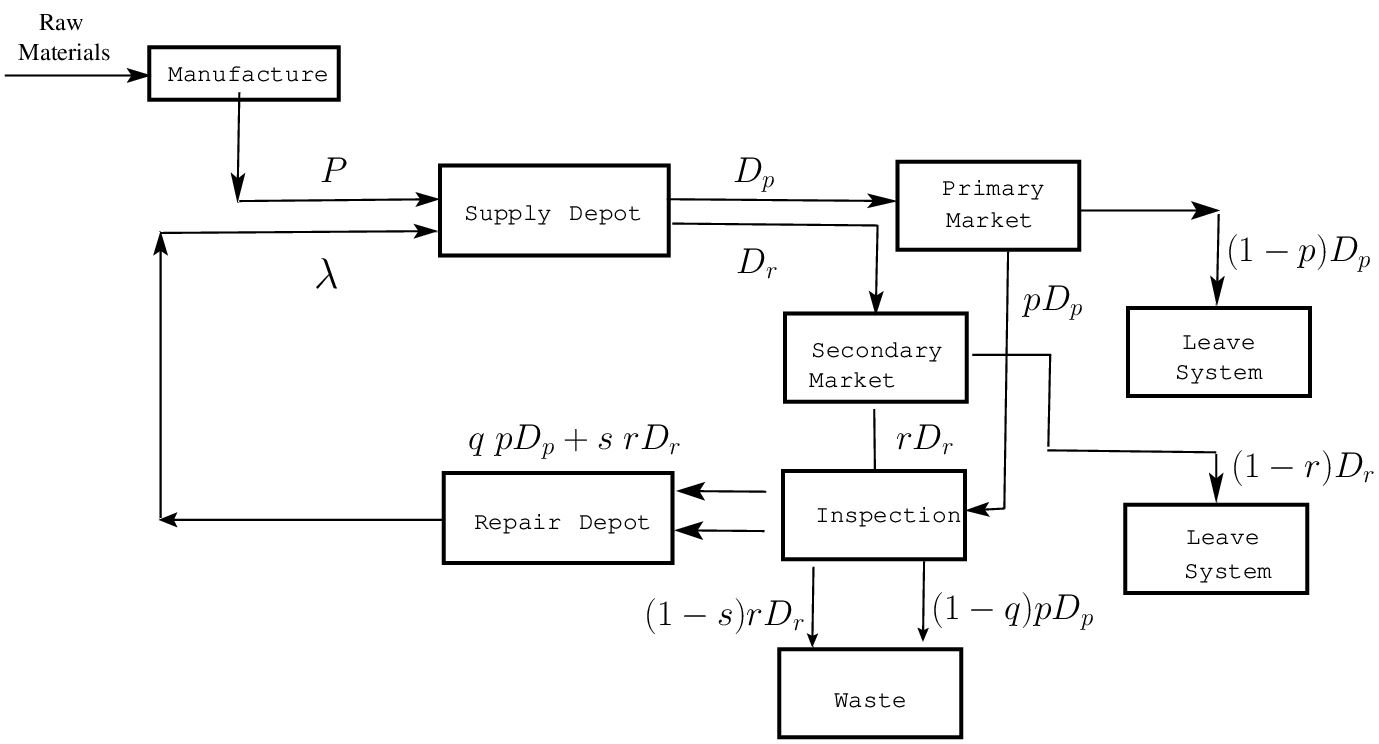}
	\caption{Material flow for a production and repair system.}
	\label{figflowModel1}
\end{figure}

The manufacturing and remanufacturing process in a reverse logistics system, as introduced by Jaber and El Saadany \cite{Jaber2009}, El Saadany \cite{Saadany2009}, Hasanov et al. \cite{Hasanov2012} and Singh and Sharma \cite{Singh2016} is illustrated in Figure \ref{figflowModel1}. A similar model was discussed by Nahmias and Rivera \cite{Nahmias1979} and Richter \cite{Richter1996a, Richter1996b, Richter1997}, with the key difference being that in \cite{Singh2016, Jaber2009, Hasanov2012} remanufactured items are sold at a discounted price in a secondary market. The system consists of two main storage areas: a supply stock, where both newly produced and remanufactured (repaired) items are stored, and a repair depot, where used items are collected. Before entering the repairable stock, used items undergo an inspection process, and non-repairable ones are discarded. Richter \cite{Richter1996b} assumed that within a given time interval \(T\), there are \(m\) remanufacturing cycles and \(n\) production cycles. Nahmias and Rivera \cite{Nahmias1979} extended the model developed by Schrady \cite{Schrady1967} by incorporating a finite repair rate at the repair depot. Their analysis also takes into account the constraints of limited storage space in both the repair and supply depots. 

\noindent In our analysis, we assume a finite repair rate \(0 <\lambda < \infty\), with the supply depot receiving repaired products at this specified rate.  Manufacturing takes place over \(n\) cycles, producing a total batch size of \(nQ_p\) during \(T_p\), while remanufacturing runs for \(m\) cycles with a total batch size of \(mQ_r\) during \(T_r\). These time intervals are defined as \(T_p = nQ_p / D_p\) and \(T_r = mQ_r / D_r\).  

\noindent In Figure \ref{figflowModel1}, the material flow system shows that items are produced by manufacturing companies at a rate $P$ and moved to the supply depot to be sold in the primary market according to customer demand $D_p$. A proportion $pD_p$ of these items are returned by customers due to faults or damage and are sent to the repair depot for sorting, while the remaining $(1-p)D_p$ items are lost from the system. After inspection, a proportion $qpD_p$ is selected for repair, while $(1-q)pD_p$ is moved to waste. Similarly, a proportion $r$ of the remanufactured or repaired items are returned, and after sorting, $srD_r$ items are selected for repair. At the repair depot, items are repaired at a rate $\lambda$. Once repaired, the items are sent to the supply depot to be sold in the secondary market at a rate $D_r$.

\subsection{Inventory cost}

\noindent In our analysis, new production and remanufacturing operations are carried out simultaneously, allowing the demand for both newly produced and repaired items to be met at the same time.  In the repair depot, returned items are collected at a combined rate of \( pqD_p + srD_r\), as returns originate from both newly sold and previously repaired items. This collection process is shown in Figure \ref{figflowModel1}.

\begin{figure}[H]
	\centering
\includegraphics[scale=1.0]{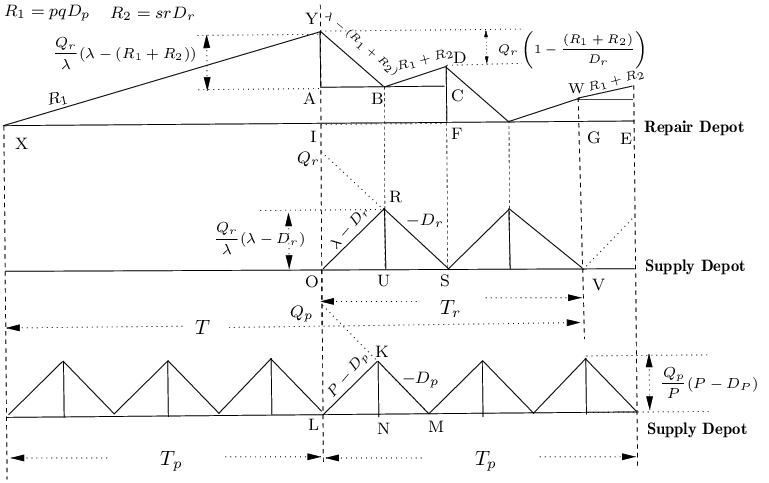}
	\caption{Inventory levels for both new and remanufactured items during the period $T_p$ and $T_r$.}
	\label{figinventory1}
\end{figure}

\noindent To evaluate the inventory accumulation, we calculate the total area under the inventory curve in the supply depot for both newly produced and remanufactured items over the periods \( T_p \) and \( T_r \), respectively. At the same time, we assess the area representing inventory buildup in the repair depot, where returned items are stored for the remanufacturing process. An important consideration here is the relative lengths of \( T_p \) and \( T_r \), as this relationship significantly affects area calculations in the repair depot. To address this complexity and simplify the analysis, we separately show the time periods and associated rates as depicted in Figure \ref{figinventory3a} (in the Appendix) and Figure \ref{figinventory2}. This allows for a more straightforward calculation while calculating the same total area as shown in Figure \ref{figinventory1}.

\noindent  In our model, we set \(\lambda \geq D_r\). In Figure \ref{figinventory2}, we observe that in the production period $T_p$, the inventory at the supply depot increases at a rate of $P - D_p$ and decreases at a rate of $D_p$. In this period, some items are returned to the repair depot, with a proportion of $pD_p$. After sorting,  $pqD_p$ items are deemed suitable for remanufacturing, resulting in an inventory increase at the repair depot at a rate of $pqD_p$ for period \(T_p\).  In the $OU$ phase, the inventory of repaired items at the supply depot increases at a rate of $\lambda - D_r$, while the inventory at the repair depot decreases at a rate of $\lambda - srD_r$, given that $\lambda \geq srD_r$. In the subsequent $US$ phase, the inventory at the supply depot decreases at a rate of $D_r$, while the inventory at the repair depot increases at a rate of $srD_r$. When the inventory level falls below \(\ds \frac{Q_r}{\lambda}(\lambda - srD_r)\) threshold, repair inductions are suspended. The net loss of repaired items over \( m \) repair cycles be \( mQ_r(1 - sr) \). The following constraints ensure that sufficient inventory is available to meet the repair demand during the period $T_r$.  
%
\begin{equation}\label{con1}
    mQ_r(1 - sr) + \frac{Q_r}{\lambda} (\lambda - D_r)sr \leq YI = npqQ_p.
\end{equation}

\begin{figure}[H]
	\centering
\includegraphics[scale=1.1]{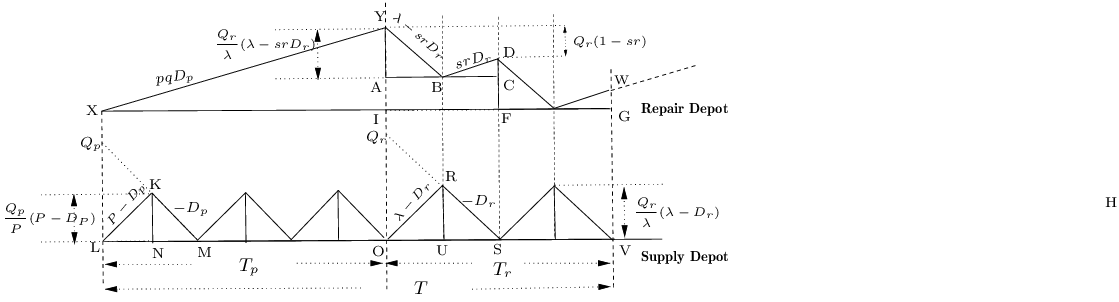}
	\caption{The behaviour of inventory for new, used and repaired items over the interval $T_p$ and $T_r$.}
	\label{figinventory2}
\end{figure}

\noindent We now derive the expression for the average area per unit time. The total cost per unit of time for the proposed reverse logistics system, denoted as \(C(Q_p, Q_r, m, n, s)\), is comprised of three main components: the setup cost per unit of time for the system, the total inventory holding cost per unit of time in the supply depot, and the total inventory holding cost per unit of time in the repair depot. We define this relation as follows. 
\begin{multline}\label{eqcost1a}
C(Q_p, Q_r, m, n, s) = \frac{1}{T}\Bigg[ mS_n+nS_p + \left(\frac{T_pQ_p(P-D_p)}{2P}+\frac{T_rQ_r(\lambda-D_r)}{2 \lambda}\right)h_p \\ +\Bigg[ \frac{R_1T_p^2}{2} + \frac{m Q_r^2(\lambda - R_2)}{2\lambda^2} +\frac{srQ_rT_r(1-\frac{D_r}{\lambda})^2}{2} \\ 
+\left(1-\frac{1}{m}\right)T_r\left(T_pR_1-\frac{Q_r}{\lambda}(\lambda - R_2)-\frac{(m-2)Q_r(1-sr)}{2}\right)\Bigg] h_r \Bigg],
\end{multline}
where, $T=T_p+T_r$, $R_1 = pqD_p$ and $R_2 = srD_r$.\\
For detailed information on area calculation, check the Appendix section \ref{appen}.

\noindent The mathematical formulation of the model is therefore written from \eqref{eqcost1a} as

\begin{equation}\label{eqmathemodel1}
\begin{array}{rl} 
\min & \ C(Q_p, Q_r, m, n, s)\\
\mbox{subject to the constraints} \\
& \ds mQ_r(1 - sr) + \frac{Q_r}{\lambda} (\lambda - D_r)sr \leq  npqQ_p, \\ 
&\ds Q_p, Q_r, m, n \in {\mathbb{Z} \geq 1},\\
& \ds 0 \leq s \leq 1.
\end{array}
\end{equation}

\subsection{Theoretical Properties of the Proposed MINLP Model}

In this section, we establish fundamental analytical properties of the proposed
mixed-integer nonlinear programming (MINLP) formulation \eqref{eqmathemodel1}, including feasibility,
boundedness, and nonconvexity of the problem.

\medskip

\noindent
Let the feasible set of problem (3) be denoted by
\[
X := \left\{ (Q_p,Q_r,m,n,s) \in \mathbb{Z}_{\ge1}^4 \times [0,1] :
mQ_r(1-sr) + \frac{Q_r}{\lambda}(\lambda-D_r)sr \le npqQ_p \right\}.
\]

We assume throughout that
\[
D_p,D_r,S_p,S_r,h_p,h_r > 0, \quad
0 < p,q,r \le 1, \quad \lambda \ge D_r.
\]


\begin{theorem}[Feasibility and Boundedness]\label{theorem1}
Under the above assumptions, the feasible set $X$ is nonempty and bounded.
\end{theorem}

\begin{proof}
\textbf{(Non-emptiness).}
Fix $m=n=1$ and choose any $s\in[0,1]$.
Constraint \eqref{con1} becomes
\[
Q_r(1-sr) + \frac{Q_r}{\lambda}(\lambda-D_r)sr \le pqQ_p.
\]
The left-hand side is linear in $Q_r$, while the right-hand side
is linear in $Q_p$. For any fixed $Q_r\ge1$, choosing $Q_p$
sufficiently large guarantees feasibility. Hence $X\neq\emptyset$.

\medskip

\textbf{(Boundedness).}
Consider the cost function $C(Q_p,Q_r,m,n,s)$ defined in \eqref{eqcost1a}.
The holding cost component contains quadratic growth terms in $Q_p$
and $Q_r$ of the form
\[
\frac{T_pQ_p(P-D_p)}{2P}, \qquad
\frac{T_rQ_r(\lambda-D_r)}{2\lambda},
\]
where $T_p = \frac{nQ_p}{D_p}$ and $T_r = \frac{mQ_r}{D_r}$.

Hence, asymptotically,
\[
C(Q_p,Q_r,m,n,s)
= \Theta(Q_p^2) + \Theta(Q_r^2) + \text{lower-order terms}.
\]
Therefore,
\[
\lim_{\|(Q_p,Q_r,m,n)\|\to\infty} C(Q_p,Q_r,m,n,s) = +\infty.
\]
Since the objective diverges to infinity as decision variables
grow unbounded, any minimizing sequence must lie in a bounded region.
Thus the feasible region containing optimal solutions is bounded.
\end{proof}

\subsubsection*{Nonconvexity of the Objective Function}

\begin{proposition}[Nonconvexity]
The objective function $C(Q_p,Q_r,m,n,s)$ is nonconvex over $X$.
\end{proposition}

\begin{proof}
Consider the production rate
\[
P(Q_p,n)
= \frac{D_p}{1 - \frac{2S_pD_p}{h_p(nQ_p)^2}}.
\]
Define
\[
\phi(Q_p) := 1 - \frac{c}{(nQ_p)^2},
\quad \text{where } c=\frac{2S_pD_p}{h_p}>0.
\]
Then
\[
P(Q_p,n) = \frac{D_p}{\phi(Q_p)}.
\]

\noindent The second derivative of $P(Q_p,n)$ with respect to $Q_p$
involves terms of order $(nQ_p)^{-4}$ and $(nQ_p)^{-6}$.
A direct computation shows that the Hessian of $C$
changes sign over the feasible domain due to the rational
structure of $P$ and $\lambda$. Hence $C$ is neither convex nor concave.
\end{proof}




\subsection{Numerical Experiments}\label{Moo}

\begin{figure}[H]
\hspace{-1cm}
\begin{center}
\includegraphics[width=110mm]{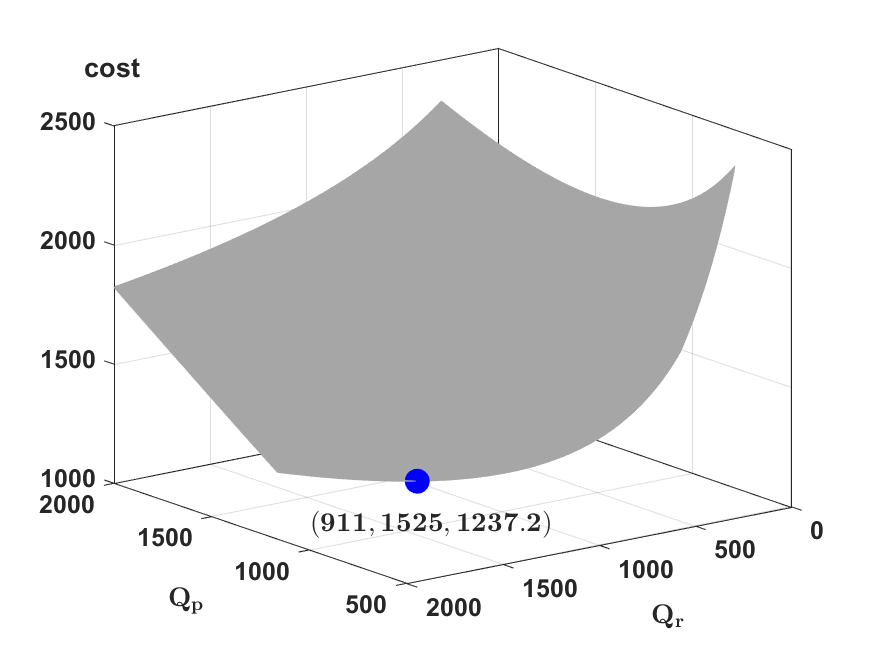} \\
\end{center}\hspace{-1cm}
\vspace{-1cm}
\caption{\small { Cost surface corresponding to Example \ref{SurfaceEPQ}, with the optimal solution highlighted at the point \( (Q_p, Q_r, \text{cost}) \).}}
\label{Ex44figsurf}
\end{figure}

\begin{example}\label{SurfaceEPQ}
This numerical illustration is based on the example provided in Sharma et al. \cite{Swati2021}. Consider a system where the primary and secondary market demands over a cycle $T_p$ and $T_r$. are both set to \( D_p = 250 \) and \( D_r = 250 \), respectively. Returned items from the primary market are collected at a rate of \( p = 0.8 \), and following a sorting process, a proportion \( q = 0.9 \) of those items are deemed recoverable. Similarly, for the secondary market, the collection rate is \( r = 0.8 \), and the post-sorting recovery efficiency is denoted by a variable \( s \), constrained within \( 0 \leq s \leq 1 \), depending on the quantity needed to fulfil the secondary market demand.

\noindent The holding costs per unit time are specified as \( h_p = 5 \) for the primary (new item) inventory and \( h_r = 2 \) for the secondary (repaired item) inventory. Setup costs are given as \( S_p = 2400 \) for new items and \( S_r = 1400 \) for repaired ones. The objective is to identify the optimal production batch sizes \( Q_p \) and \( Q_r \), as well as the number of cycles \( n \) for new items and \( m \) for repaired items, that together minimize the average inventory cost over the cycle time $T_p$ and $T_r$.. This is achieved by minimizing the cost function defined in \(C(Q_p, Q_r, m, n, s)\), subject to the constraints outlined in \eqref{eqmathemodel1}. Figure \ref{Ex44figsurf} illustrates the surface generated by \eqref{eqmathemodel1}, with the aim of demonstrating its optimal solution under the given conditions.
\end{example}

\noindent Based on the parameter settings in Example \ref{SurfaceEPQ}, the model yields a minimum total cost of \( 1237.2 \), with optimal batch sizes of \( Q_p = 911 \) for new items and \( Q_r = 1525 \) for repaired items. The corresponding cycle frequencies are \( n = 1 \) and \( m = 2 \), and a full recovery rate \( s = 1 \) is required to meet the secondary market demand.

\subsection{Sensitivity Analysis} \label{sec4b}
\begin{figure}[H]
\hspace{-1cm}
\begin{minipage}{100mm}
\begin{center}
\hspace*{-2.5cm}
\includegraphics[width=80mm]{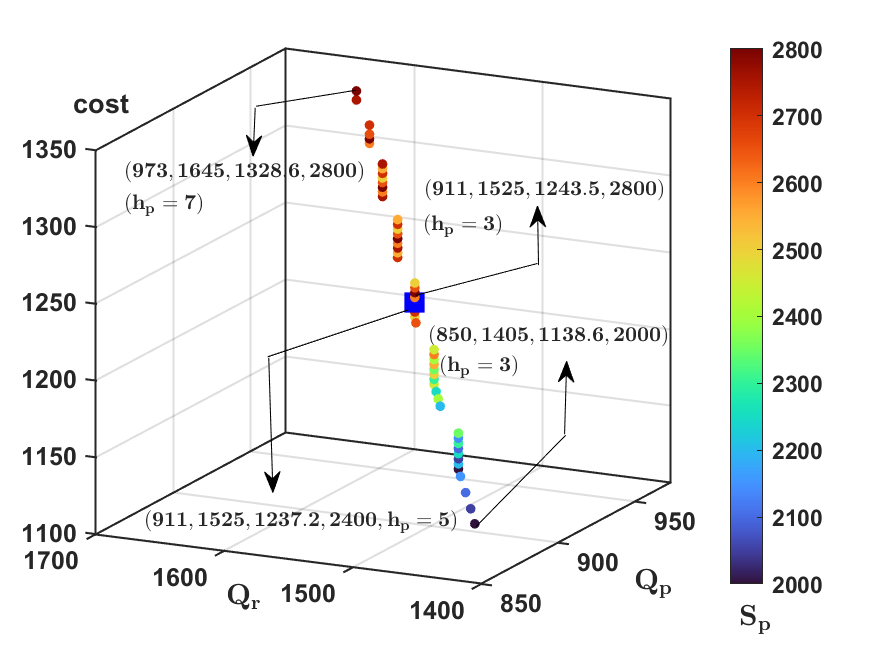} \\
\hspace*{-2.9cm}
{\scriptsize (a) Inventory cost with respect to $Q_p$, $Q_r$, and $S_p$.}
\end{center}
\end{minipage}
\hspace*{-3.5cm}
\begin{minipage}{100mm}
\begin{center}
\hspace*{0cm}
\includegraphics[width=80mm]{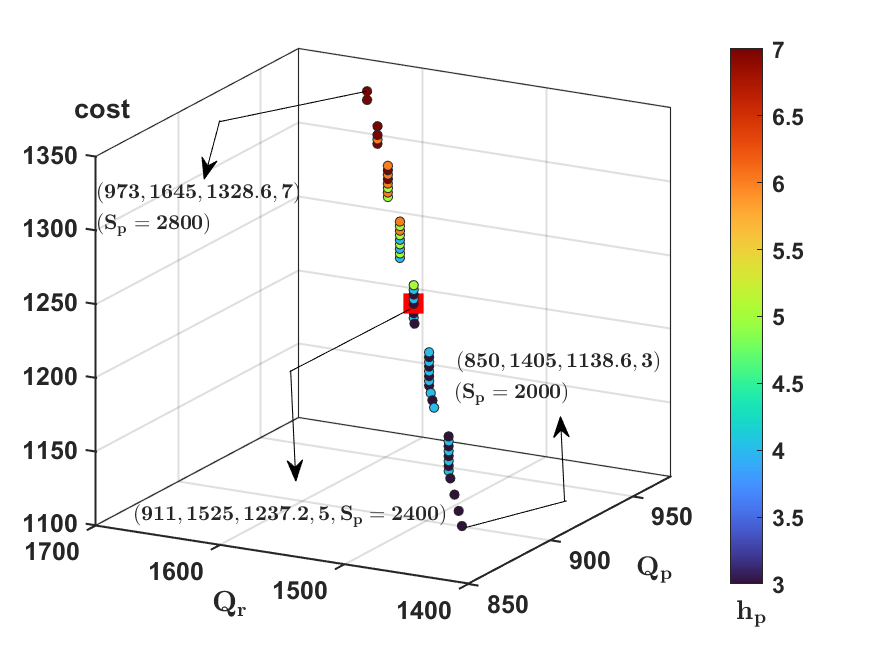} \\
\hspace*{-1.0cm}
{\scriptsize (b) Inventory cost with respect to $Q_p$, $Q_r$, and $h_p$.}
\end{center}
\end{minipage}
\caption{{Inventory cost across a range of setup $(S_p)$ and holding $(h_p)$ costs.
}}
\label{figSphp}
\end{figure}

\begin{figure}[H]
\hspace{-1cm}
\begin{minipage}{100mm}
\begin{center}
\hspace*{-2.5cm}
\includegraphics[width=80mm]{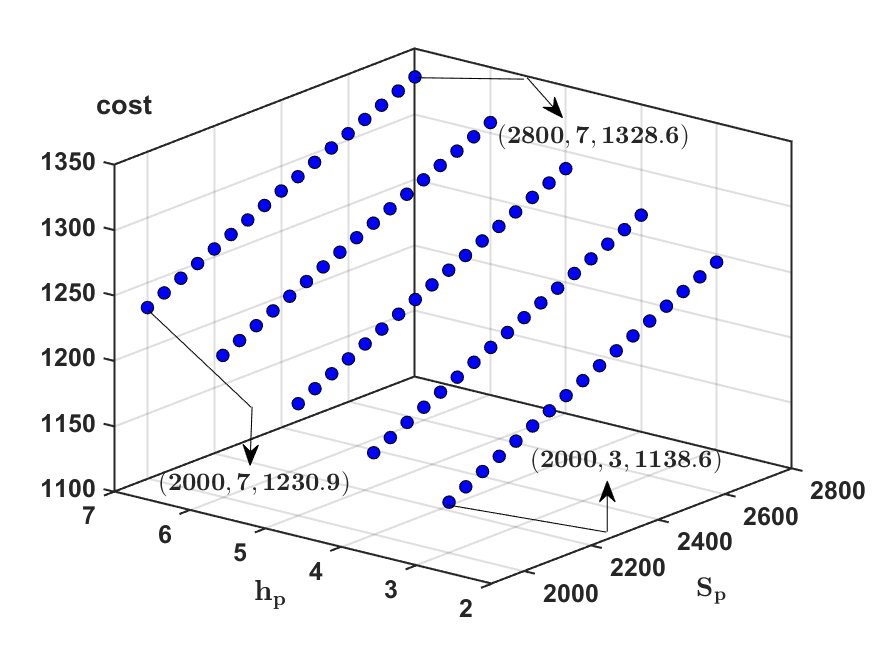} \\
\hspace*{-2.9cm}
{\scriptsize (a)  Inventory cost with respect to $h_p$ and $S_p$.}
\end{center}
\end{minipage}
\hspace*{-3.5cm}
\begin{minipage}{100mm}
\begin{center}
\hspace*{0cm}
\includegraphics[width=80mm]{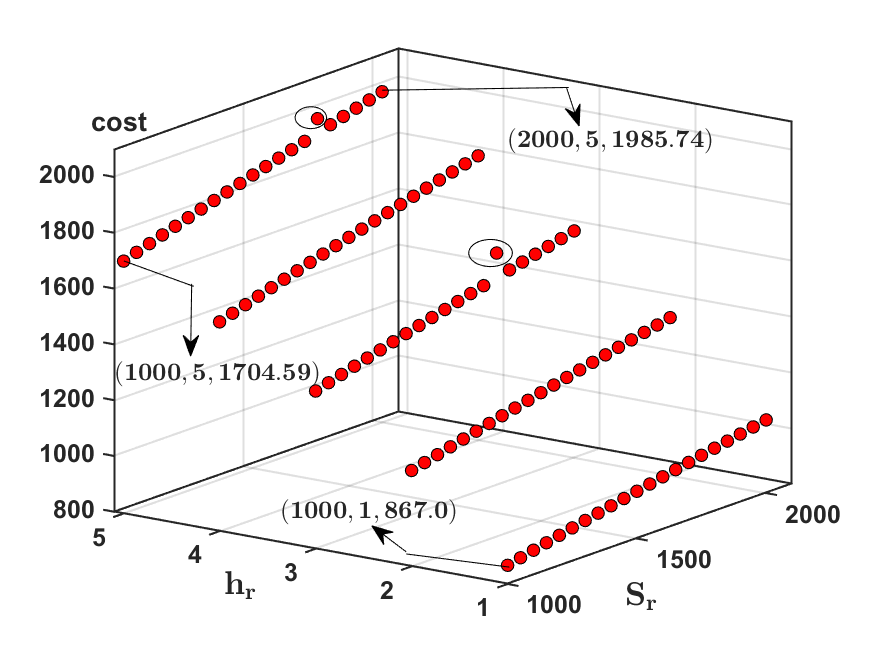} \\
\hspace*{-1.0cm}
{\scriptsize (b) Inventory cost with respect to $h_r$ and $S_r$.}
\end{center}
\end{minipage}
\caption{{Effect of setup and holding costs on inventory cost: comparing $(S_p, h_p)$ and $(S_r, h_r)$ scenarios.
}}
\label{fig_setup_holding_cost}
\end{figure}

\begin{figure}[H]
\hspace{-1cm}
\begin{minipage}{100mm}
\begin{center}
\hspace*{-2.5cm}
\includegraphics[width=80mm]{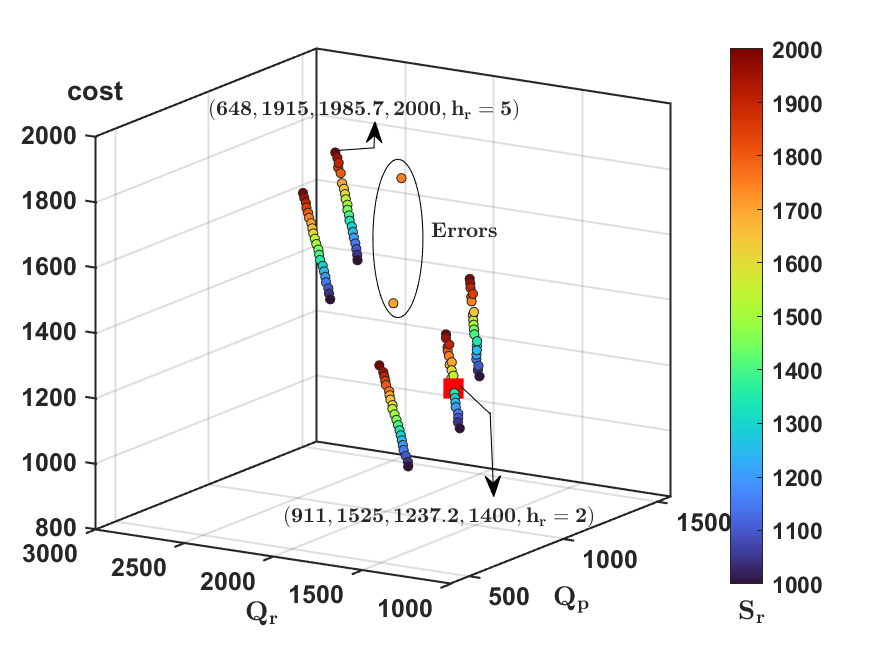} \\
\hspace*{-2.9cm}
{\scriptsize (a) Inventory cost with respect to $Q_p$, $Q_r$, and $S_r$.}
\end{center}
\end{minipage}
\hspace*{-3.5cm}
\begin{minipage}{100mm}
\begin{center}
\hspace*{0cm}
\includegraphics[width=80mm]{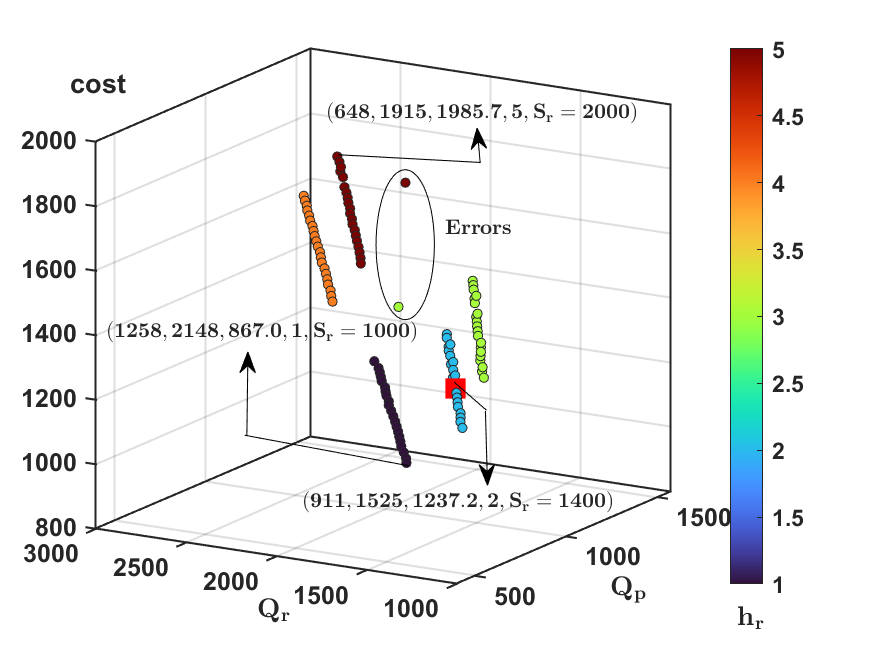} \\
\hspace*{-1.0cm}
{\scriptsize (b) Inventory cost with respect to $Q_p$, $Q_r$, and $h_r$.}
\end{center}
\end{minipage}
\caption{{Inventory cost across a range of setup $(S_r)$ and holding $(h_r)$ costs.
}}
\label{fighpCost}
\end{figure}

A sensitivity analysis is carried out to evaluate the effects of changes in setup and holding costs for both production and remanufacturing on the optimal solution. The results are depicted in Figures \ref{figSphp} to \ref{fighpCost}. The key findings are summarised below.

\begin{enumerate}
\item [(i)] Figure~\ref{figSphp}(a) and (b) illustrate the variation in production setup cost \(S_p\) from 2000 to 2800 and supply depot holding cost \(h_p\) from 3 to 7, while keeping all other parameters constant as specified in Example 5.1. In Figure~\ref{figSphp}(a), when \((S_p, h_p) = (2000, 3)\), the optimal lot sizes for production and remanufacturing are 850 and 1405 items, respectively, with an optimal inventory cost of \(\$1138.6\). However, as \(S_p\) increases to \(\$2800\) (\(h_p=3\)) per order, the optimal lot sizes rise to \(911\) items for production and \(1525\) items for remanufacturing, resulting in an increased inventory cost of \(\$1243.5\). This analysis reveals that while the setup cost increases by $40\%$, the inventory cost rises by approximately $9.23\%$. Based on the results, we cannot conclude that increasing the lot size will necessarily lead to minimum cost, as other factors—such as setup costs and the holding costs of repaired items, total cycle time—also influence overall inventory costs. In addition, we also see that if the holding cost increases by $40\%$, the inventory cost is estimated to increase by approximately $2.43\%$.

\item [(ii)] In Example \ref{SurfaceEPQ}, we consider a setup cost \( S_p = 2400 \) and holding cost \( h_p = 5 \). This setting results in an associated cost of $\$ 1237.2$, as indicated by the square marker in Figure \ref{figSphp}(b). From Figures \ref{figSphp}(b) and \ref{fig_setup_holding_cost}(a), we observe how different combinations of $(S_p, h_p)$ influence the total cost. The example illustrates the importance of selecting an appropriate balance between setup and holding costs in order to minimize overall inventory cost. This balance is crucial in production planning, as higher setup costs tend to favour larger, less frequent production runs, while higher holding costs encourage smaller, more frequent runs. Therefore, the combination shown in Example 5.1 serves as a reference for making cost-efficient decisions in inventory management.
\item [(iii)] We performed a similar analysis for the setup and holding costs related to repaired or remanufactured products, as we showed for $S_p$ and $h_p$. In this case, the setup cost $S_r$ varies from $1000$ to $2000$, and the holding cost $h_r$ ranges from $1$ to $5$, while all other parameters are kept the same. The results are shown in Figures \ref{fig_setup_holding_cost}(b), \ref{fighpCost}(a), and (b). In Example \ref{SurfaceEPQ}, we focus on the combination $(S_r, h_r) = (1400, 2)$, which is marked by a square in Figures \ref{fighpCost}(a) and (b). The optimal results for other combinations of $S_r$ and $h_r$ are also depicted in these figures. We solved the model for each pair, and the solutions are shown in the figures. We observe that the improved results are obtained for $(S_r, h_r) = (1000, 1)$, highlighting the importance of choosing the appropriate combination between setup and holding costs to minimize the total inventory cost.
 Some areas in the figures are outlined with elliptical shapes. These indicate errors where the solver could not find the optimal solution. This happens because of a local solver, which can sometimes get stuck at a local minimum instead of finding the global solution. A variety of initial guesses can help converge on the solution.

\end{enumerate}

\subsection{Real-World Applicability}
To assess the practical relevance of the proposed model, we analyzed its behavior under parameter ranges consistent with real industrial settings such as automotive components, refurbished electronics, and tyre retreading industries.
  In such environments, firms simultaneously manage new product manufacturing and remanufacturing operations under quality-differentiated demand streams. Using industry-informed parameter ranges, the model was aligned to reflect realistic operational conditions. The results confirm that incorporating integer lot-sizing decisions and variable production rates yields implementable batch policies aligned with practical production planning constraints. Furthermore, the multiobjective framework provides decision-makers with explicit trade-off curves between cost efficiency and environmental impact, supporting compliance with sustainability regulations and corporate ESG targets. These findings demonstrate that the proposed mixed-integer reverse logistics EPQ model is not only theoretically sound but also operationally viable for manufacturing firms transitioning toward circular economy practices.

\section{Environmental Impact}\label{EnvImp}
In this section, we define the greenhouse gas (GHG) emissions and energy consumption associated with both production and remanufacturing processes. These environmental impacts are functions of the production and remanufacturing rates, which vary based on lot sizes and cycle lengths. Later, in Section \ref{MOP}, we will combine these emission and energy functions with the inventory cost model to study the trade-off between cost and environmental impact.

\subsection{GHG emission}
Jaber et al. \cite{Jaber2013a} and Bazan et al. \cite{Bazan2015, Bazan2016} considered the amount of green house gas emissions from production and remanufacturing activities as $G_p(P)=a_p P^2 - b_p P + c_p$ and $G_r(\lambda)=a_r \lambda^2 - b_r \lambda + c_r$, respectively.

\noindent Substituting the production and remanufacturing rate taken from \cite{Bazan2016}
\(\ds P(Q_p, n) = \frac{D_p}{1 - \frac{2S_pD_p}{h_p(nQ_p)^2}}\)
and \(\ds \lambda(Q_r, m) = \frac{D_r}{1 - \frac{2S_rD_r}{h_r(mQ_r)^2}}\), we have
$$
G_p(Q_p, n) = a_p \left(\frac{D_p}{1 - \frac{2S_pD_p}{h_p(nQ_p)^2}}\right)^2 - b_p \left(\frac{D_p}{1 - \frac{2S_pD_p}{h_p(nQ_p)^2}}\right) + c_p,
$$
and

$$
G_r(Q_r, m) = a_r \left(\frac{D_r}{1 - \frac{2S_rD_r}{h_r(mQ_r)^2}}\right)^2 - b_r \left(\frac{D_r}{1 - \frac{2S_rD_r}{h_r(mQ_r)^2}}\right) + c_r.
$$
From Bazan et al. \cite{Bazan2015}, the annual GHG emissions cost from transportation is estimated based on the emissions generated per gallon of diesel fuel consumed. In this study, transportation involves trucks running on diesel with an average fuel efficiency of $4$ miles per gallon (mpg), travelling a round-trip distance of $300$ miles between the production facility and the market. Given this setup, the total gallons of fuel consumed per trip is $\ds \frac{300}{4} = 75$ gallons. The annual transportation-related GHG emissions cost for both production and remanufacturing is estimated using the following expression:

$$
G_t = \left( \frac{D_p + D_r}{t_c} \right) \cdot g_t \cdot e_t,
$$

\noindent where $D_p$ and $D_r$ are the annual product and remanufactured product demands, respectively, $t_c$ is the number of units transported per truck trip (truck capacity), $g_t$ is the gallons of diesel consumed per trip ($75$ gallons in this case), and $e_t$ is the GHG emissions cost per gallon of diesel fuel. This formulation captures the emissions cost based on the number of trips required to meet total annual demand and the emissions produced per trip. Therefore, the total cost associated with GHG emissions from production, remanufacturing, and the transportation related to both processes is as follows.

\begin{equation}
    G(Q_p, Q_r, m , n)= \left(G_p(Q_p, n)+G_r(Q_r, m)+G_t \right)C_{GHG}.
\end{equation}


\begin{proposition}[Lot Size–Emission Monotonicity]
Assume the emission function
\[
G_p(P) = a_pP^2 - b_pP + c_p,
\quad a_p>0,
\]
is convex in $P$. Then for sufficiently large $Q_p$,
$G_p(Q_p,n)$ is decreasing in $Q_p$.
\end{proposition}

\begin{proof}
Since
\[
P(Q_p,n)
= \frac{D_p}{1 - \frac{c}{(nQ_p)^2}},
\]
we compute
\[
\frac{\partial P}{\partial Q_p}
= -\frac{D_p \cdot 2c}{(n^2Q_p^3)\left(1 - \frac{c}{(nQ_p)^2}\right)^2} < 0.
\]
Thus $P$ is strictly decreasing in $Q_p$.

By the chain rule,
\[
\frac{\partial G_p}{\partial Q_p}
= (2a_pP - b_p)\frac{\partial P}{\partial Q_p}.
\]
For sufficiently large $Q_p$, $P$ remains bounded,
and $(2a_pP - b_p)$ is positive under standard parameter settings.
Since $\frac{\partial P}{\partial Q_p}<0$,
we obtain
\[
\frac{\partial G_p}{\partial Q_p} < 0.
\]
Hence emissions decrease as lot size increases.
\end{proof}

\medskip

\noindent This formally establishes the structural trade-off:
larger lot sizes reduce emissions through endogenous rate adjustment,
while increasing holding costs.

\subsection{Energy used for production}

Zanoni et al. \cite{Zanoni2014} and later Bazan et al. \cite{Bazan2015, Bazan2016} and Forkan et al. \cite{Forkan2022} introduced mathematical formulations to estimate the annual energy consumption in both production and remanufacturing processes. The energy used in the production process, denoted as 
$$
E_p(Q_p, n) = \left(K_p + \frac{2400W_p\left(1 - \frac{2S_pD_p}{h_p(nQ_p)^2}\right)}{D_p}\right)nQ_p,
$$ 
and in the remanufacturing process, denoted as 
$$
E_r(Q_r, m) = \left(K_r + \frac{2400W_r\left(1 - \frac{2S_rD_r}{h_r(mQ_r)^2}\right)}{D_r}\right)mQ_r,
$$
are calculated using expressions that account for both machine-specific energy requirements and operational parameters. Specifically, the equations incorporate the fixed energy required to produce a single unit $K_p$ (kWh/unit) and $K_r$ (kWh/unit) for production and repair production, respectively. The idle power consumption of the machines are $W_p$ (kW) and $W_r$ (kW).  These formulas are based on an annual operation of $2400$ hours, assuming $300$ working days at $8$ hours per day. 

\noindent Therefore, the total cost associated with energy used in production and remanufacturing is as follows.

\begin{equation}
    E(Q_p, Q_r, m , n)= \left(E_p(Q_p, n)+E_r(Q_r, m) \right)C_{ENG}.
\end{equation}


%


\section{Multiobjective Optimization Problems} \label{MOP}
\noindent Multiobjective optimization problems (MOPs) involve the simultaneous optimization of multiple, often conflicting, objectives. Due to the inherent trade-offs between these objectives, it's typically impossible to find a single solution that optimizes all objectives simultaneously. Instead, the goal is to identify a set of optimal solutions, known as the Pareto front, where no objective can be improved without compromising at least one other. A common strategy for addressing MOPs is scalarization, which transforms a multiobjective problem into a single objective one. This transformation facilitates the use of traditional optimization techniques. One such scalarization method is the feasible value constraint approach \cite{BurKayRiz2017}, which focuses on optimizing a specific objective while treating the remaining objectives as constraints. This method allows for the approximation of Pareto-optimal solutions by systematically varying the constraints imposed on the feasible set defined in \eqref{eqmathemodel1}. More approaches for scalarization can also be found in the literature (see, \cite{PascolettiSerafini1984}--\cite{Burachik2014}). In our analysis, we employed the feasible value constraint approach due to its straightforward implementation and effectiveness in generating Pareto-optimal solutions. 

\noindent The objective of this study is not to develop a new scalarization methodology but rather to demonstrate how established multiobjective optimization techniques can be effectively integrated with reverse logistics inventory systems involving endogenous production and remanufacturing rates. The feasible-value-constraint approach was selected because of its ability to generate Pareto-optimal solutions for nonconvex optimization problems and mixed-integer formulations while maintaining computational tractability.

\noindent The descriptions of the approach for two objectives are given below.

\subsection{Multiobjective formulation with inventory and GHG emissions}
\noindent 
Most existing reverse logistics inventory models focus primarily on cost minimization, treating environmental impacts as secondary considerations or evaluating them separately after optimization. In contrast, the proposed framework explicitly incorporates environmental objectives into the decision-making process, allowing managers to evaluate trade-offs between economic and environmental performance simultaneously. Such trade-off analysis is increasingly important in manufacturing systems operating under sustainability targets, ESG reporting requirements, and circular economy initiatives.

\noindent We begin by formulating a two-objective optimization model that simultaneously addresses two critical aspects of production systems: inventory-related costs and greenhouse gas (GHG) emissions arising from both manufacturing and remanufacturing processes. This approach aims to identify optimal trade-offs between economic efficiency and environmental sustainability, and formulate as follows:

\begin{equation}\label{modelmop1}
\begin{array}{rl} 
\min & \ \bigg ( C(Q_p, Q_r, m, n, s), G(Q_p, Q_r, m, n) \bigg )\\
\mbox{subject to the constraints} \\
& \ds mQ_r(1 - sr) + \frac{Q_r}{\lambda} (\lambda - D_r)sr \leq  npqQ_p, \\ 
&\ds Q_p, Q_r, m, n \in {\mathbb{Z} \geq 1},\\
& \ds 0 \leq s \leq 1.
\end{array}
\end{equation}

\noindent We utilize \textbf{\em the feasible-value-constraint approach \em\,:} For $(\hat{Q}_p,\hat{Q}_r,\hat{m},\hat{n},\hat{s}) \in X$, where $X$ be a feasible set of constraints stated in \eqref{modelmop1}, and
\(\ds W := \left
\{ w \in \mathbb{R}^{2} \mid w_1, w_2 > 0,\sum_{i=1}^{2} w_i=1
\right \}.\)
\noindent By varying $w$ within this set, we can explore different trade-offs and construct the Pareto front, which represents the set of non-dominated solutions where no objective can be improved without worsening the other.

\noindent We assume $w_1C(\hat{Q}_p,\hat{Q}_r,\hat{m},\hat{n},\hat{s})=w_2G(\hat{Q}_p,\hat{Q}_r,\hat{m},\hat{n})$. We define two scalar problems as
\begin{equation} \label{nchep1}
\mbox{($P_{\hat{x}}^1$)}\ \left\{\begin{array}{rl} \ds\min_{(Q_p,Q_r,m,n,s) \in  X} & \
\ w_1C(Q_p,Q_r,m,n,s),
\\[4mm]
\mbox{subject to} & \ \ w_2G(Q_p, Q_r, m , n) \leq w_1C(\hat{Q}_p,\hat{Q}_r,\hat{m},\hat{n},\hat{s}),
\end{array}
\right.
\end{equation}
and
\begin{equation} \label{nchep2}
\mbox{($P_{\hat{x}}^2$)}\ \left\{\begin{array}{rl} \ds\min_{(Q_p,Q_r,m,n,s)\in  X} & \
\ w_{2}G(Q_p, Q_r, m , n),
\\[4mm]
\mbox{subject to} & \ \ w_1C(Q_p,Q_r,m,n,s) \leq w_2G(\hat{Q}_p, \hat{Q}_r, \hat{m} , \hat{n}).
\end{array}
\right.
\end{equation}

\noindent
\begin{example}\label{multexam}
	We proceed to solve the bi-objective optimization problems \eqref{nchep1} and \eqref{nchep2} using the following parameter settings:  \(D_p=1000\) units/month, \(D_r=422\) units/month, \(p=0.6\) ($\%$), \(q=0.9\) ($\%$), \(r=0.7\) ($\%$), \(S_p=\$50\)/setup, \(S_r=\$40\)/setup,  \ and \(h_p=h_r=\$15\)/unit/month, $a_p=0.003$ $ton \;year^2$/$unit^3$, $b_p=0.12$ $ton \;year$/$unit^2$ $c_p=1.4$ $ton$/$unit$, $a_r=0.003$ $ton \;year^2$/$unit^3$, $b_r=0.12$ $ton \;year$/$unit^2$, $c_r=1.4$ $ton$/$unit$ \ and \ $C_{GHG}=\$ 2$;
\end{example}
\noindent The objective is to generate a set of Pareto-optimal solutions that illustrate the trade-offs between inventory-related costs and greenhouse gas (GHG) emissions. To obtain this, we solve the problems introduced in \eqref{nchep1} and \eqref{nchep2} using the algorithm stated in  \cite[Appendix A.2., Algorithm 1]{Indu2025}.

\begin{figure}[H]
\hspace{-1.1cm}
\begin{minipage}{82mm}
\begin{center}
\hspace*{0cm}
\includegraphics[width=82mm]{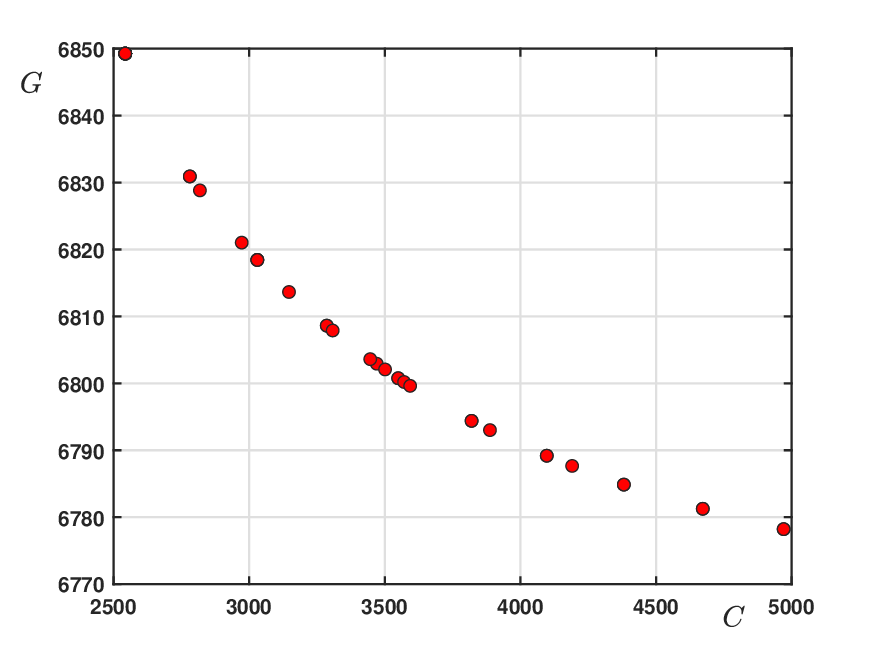} \\
{\scriptsize (a) Pareto front approximated with $P_{\hat{x}}^k$.}
\end{center}
\end{minipage}
\hspace{-0.7cm}
\begin{minipage}{82mm}
\begin{center}
\hspace*{0cm}
\includegraphics[width=82mm]{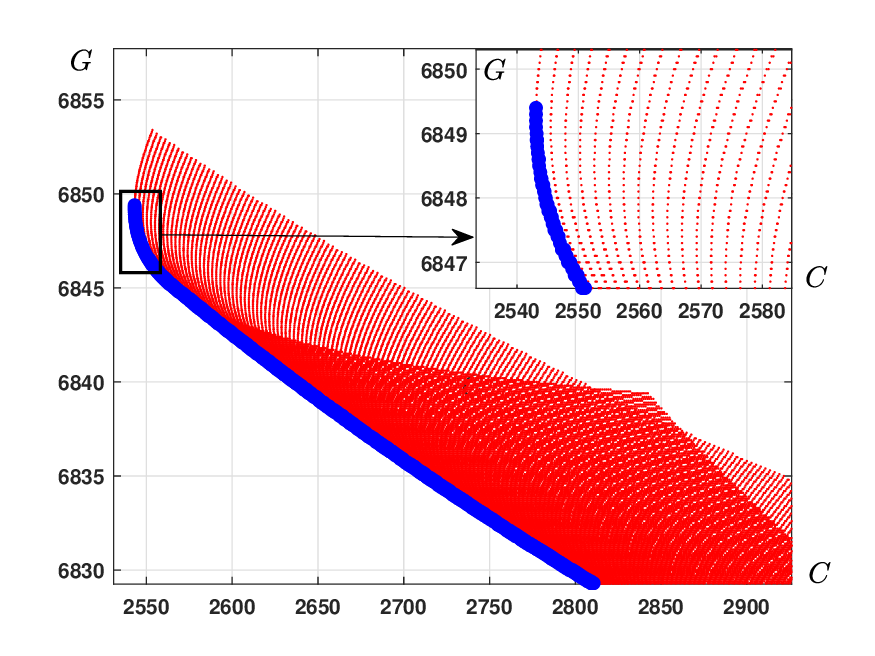} \\
{\scriptsize (b) Heuristic approach is used to construct the Pareto front.}
\end{center}
\end{minipage}
\caption{{Equal grids are provided in the Algorithm.
            (Red) circles indicate the Pareto front obtained by $P_{\hat{x}}^k$ in (a), (blue) circles depict Pareto points found by the brute-force algorithm in (b), and Pareto points are displayed on the magnified axis in (b). }}
\label{multiexampleA}
\end{figure}

\noindent In our results, we observe that the lowest inventory cost of $\$2,543.13$ is achieved when the lot sizes are set to $Q_p = 1000$ and $Q_r = 566$, with production cycles $m = n = 1$. However, this combination leads to the highest GHG cost of $\$6,849.22$. We also observe that the high GHG cost corresponds to a production rate $P = 1006$ and a remanufacturing rate $\lambda = 425$. As we increase the lot sizes $Q_p$ and $Q_r$, the production and remanufacturing rates decrease accordingly. This reduction in $P$ and $\lambda$ leads to lower GHG emissions, as the emission function $G(.)$ is convex and sensitive to high production rates. However, this improvement in environmental performance comes at the cost of higher inventory costs, since larger lot sizes result in more inventory being held over time.

\noindent This trade-off—lower GHG emissions at the expense of higher inventory cost, and vice versa—is illustrated in Figure \ref{multiexampleA}(a), highlighting the conflicting nature of these two objectives. In Figure \ref{multiexampleA}(a), we use BONMIN to solve the subproblems defined in \eqref{nchep1} and \eqref{nchep2}, following the algorithm proposed in \cite[Appendix A.2., Algorithm 1]{Indu2025} to approximate the Pareto points. These Pareto points correspond to various weight vectors $w \in W$, which must be specified as input to the algorithm. Although we tested with other mixed-integer solvers, BONMIN \cite{Bonami2008} demonstrated better performance in solving the subproblems \eqref{nchep1} and \eqref{nchep2}. Table \ref{tablemop1} presents a summary of the number of weight vectors, reference points, subproblems solved, and the resulting Pareto and non-Pareto points. Non-Pareto points appear when the solver becomes stuck at infeasible points. These non-dominated points are subsequently filtered out using elimination techniques.
\begin{table}[H]
\small
\caption{Numerical performance of $P_{\hat{x}}^k$.}  
\centering 
\begin{tabular}{||c| c| c| c| c| c||} 
\hline \hline 
Methods &   Number of  & Number of   & Time          & Pareto    & Non-Pareto  \\  
       &   grid points     & subproblems & in seconds    & points    &   points \\ [0.5ex] 
\hline 
$P_{\hat{x}}^k$ & 50       & 100       & 1342       & 52        & 18          \\ 
\hline
\hline 
\end{tabular}
\label{tablemop1} 
\end{table}

\noindent Figure \ref{multiexampleA}(b) presents the objective space defined by inventory cost and GHG emissions, with the goal of visualizing the Pareto front and verifying the results shown in Figure \ref{multiexampleA}(a). In this analysis, we consider the lot size ranges $Q_p \in [1000, 1500]$ and $Q_r \in [500, 1000]$, and evaluate both the inventory cost $C(\cdot)$ and GHG cost $G(\cdot)$ using $m = n = s = 1$. This parameter space results in $251,001$ unique combinations of $(C, G)$ values, which are illustrated as the red surface in Figure \ref{multiexampleA}(b).
To identify the Pareto-optimal solutions, we apply a brute-force algorithm to eliminate dominated solutions. This process yields $11,320$ non-dominated points, which form the Pareto front, shown as the blue curve in the figure. 
We then compare these results with Figure \ref{multiexampleA}(a), which is generated using a limited number of grids. The comparison shows that the Pareto points obtained in Figure \ref{multiexampleA}(a) are indeed a subset of the more comprehensive set of non-dominated points identified in Figure \ref{multiexampleA}(b).

\subsection{Multiobjective formulation with inventory and energy consumptions}
\noindent 
Another important aspects of production systems: inventory-related costs and energy consumption arising from both manufacturing and remanufacturing processes. The purpose of this multiobjective model is to identify optimal trade-offs between minimizing total inventory related costs and reducing energy usage, and it is formulated as follows:

\begin{equation}\label{modelmop2}
\begin{array}{rl} 
\min & \ \bigg ( C(Q_p, Q_r, m, n, s), E(Q_p, Q_r, m, n) \bigg )\\
\mbox{subject to the constraints} \\
& \ds mQ_r(1 - sr) + \frac{Q_r}{\lambda} (\lambda - D_r)sr \leq  npqQ_p, \\ 
&\ds Q_p, Q_r, m, n \in {\mathbb{Z} \geq 1},\\
& \ds 0 \leq s \leq 1.
\end{array}
\end{equation}

\begin{figure}[H]
\hspace{-1.1cm}
\begin{minipage}{82mm}
\begin{center}
\hspace*{0cm}
\includegraphics[width=82mm]{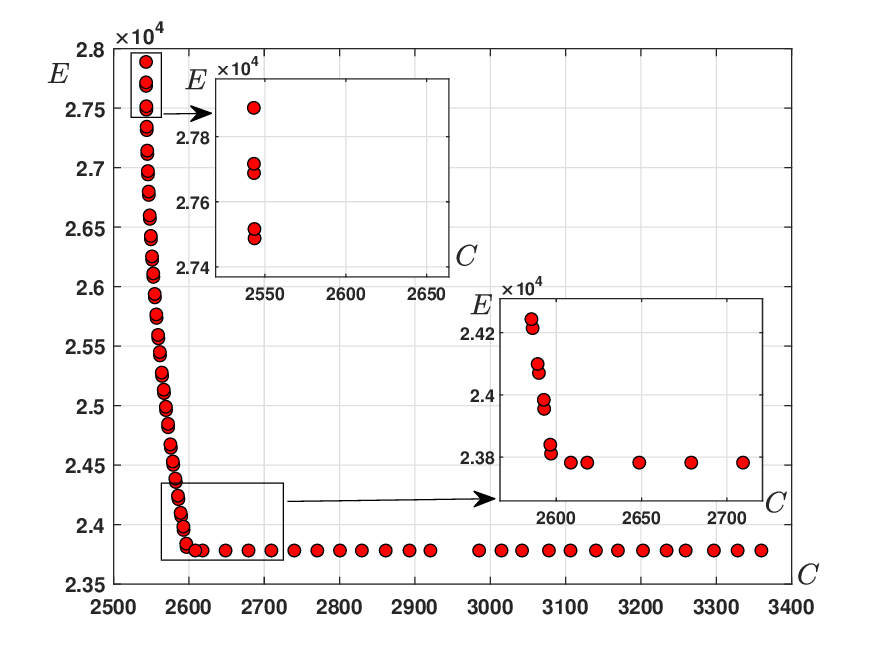} \\
{\scriptsize (a) Pareto front approximated with $P_{\hat{x}}^k$.}
\end{center}
\end{minipage}
\hspace{-0.7cm}
\begin{minipage}{82mm}
\begin{center}
\hspace*{0cm}
\includegraphics[width=82mm]{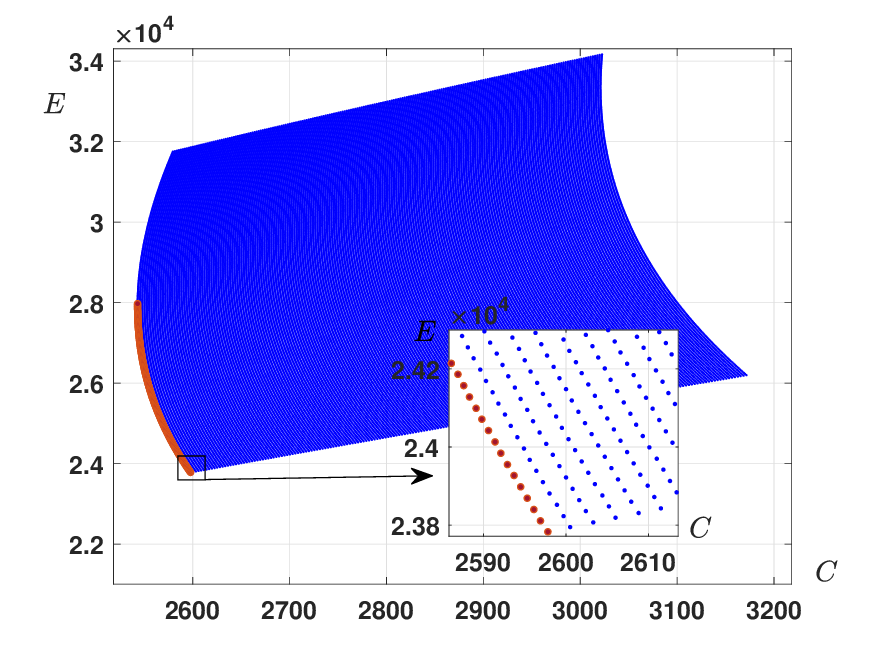} \\
{\scriptsize (b) Heuristic approach is used to approximate the Pareto front.}
\end{center}
\end{minipage}
\caption{{Equal grids are provided in the Algorithm.
            (Red) circles indicate the Pareto front obtained by $P_{\hat{x}}^k$ in (a), (blue) circles depict Pareto points found by the brute-force algorithm in (b), and Pareto points are displayed on the magnified axis in (b).  }}
\label{multiexample}
\end{figure}

\noindent \noindent We again apply \textbf{\em the feasible-value-constraint approach}, considering all $w \in W$ and assume
$w_1C(\hat{Q}_p,\hat{Q}_r,\hat{m},\hat{n},\hat{s})=w_2E(\hat{Q}_p,\hat{Q}_r,\hat{m},\hat{n})$. We define two scalar problems as
\begin{equation} \label{nchep3}
\mbox{($P_{\hat{x}}^3$)}\ \left\{\begin{array}{rl} \ds\min_{(Q_p,Q_r,m,n,s) \in  X} & \
\ w_1C(Q_p,Q_r,m,n,s),
\\[4mm]
\mbox{subject to} & \ \ w_2E(Q_p, Q_r, m , n) \leq w_1C(\hat{Q}_p,\hat{Q}_r,\hat{m},\hat{n},\hat{s}),
\end{array}
\right.
\end{equation}
and
\begin{equation} \label{nchep4}
\mbox{($P_{\hat{x}}^4$)}\ \left\{\begin{array}{rl} \ds\min_{(Q_p,Q_r,m,n,s)\in  X} & \
\ w_{2}E(Q_p, Q_r, m , n),
\\[4mm]
\mbox{subject to} & \ \ w_1C(Q_p,Q_r,m,n,s) \leq w_2E(\hat{Q}_p, \hat{Q}_r, \hat{m} , \hat{n}).
\end{array}
\right.
\end{equation}
\noindent We set the reference point (utopia vector) as $(-1000, -1000)$ and generated uniformly distributed weight vectors $w_i$ to systematically explore the trade-offs between inventory-related costs and energy consumption-related costs. This approach approximated $200$ solutions. However, approximately $20\%$ of these were identified as dominated points which are removed from the front (see, Table \ref{tablemop2}). In Figure \ref{multiexample}(a), we observe that the minimum inventory cost is \(2,543\), the corresponding energy consumption cost is \(27,516\). Conversely, when the energy cost is minimized at \(23,782\), the associated inventory cost increases to \(2,618\). The resulting Pareto front is illustrated in Figure \ref{multiexample}(a).

\begin{table}[H]
\small
\caption{Numerical performance of $P_{\hat{x}}^k$.}  
\centering 
\begin{tabular}{||c| c| c| c| c| c||} 
\hline \hline 
Methods &   Number of  & Number of   & Time          & Pareto    & Non-Pareto  \\  
       &   grid points     & subproblems & in seconds    & points    &   points \\ [0.5ex] 
\hline 
$P_{\hat{x}}^k$ & 50       & 100       & 1642       & 140        & 20          \\ 
\hline
\hline 
\end{tabular}
\label{tablemop2} 
\end{table}
\noindent The results of this study reveal a clear trade-off between inventory cost and energy consumption. Specifically, when inventory cost is minimized, energy consumption tends to be higher, and vice versa. This occurs because minimizing inventory cost often involves producing in smaller lot sizes or more frequently, which reduces holding costs but increases the number of production setups and overall energy use. Conversely, minimizing energy consumption usually requires producing in larger batches with fewer setups, which leads to higher inventory levels and increased holding costs. This inverse relationship highlights the need for a balanced approach, where decision-makers must weigh cost efficiency against environmental sustainability. The trade-off is further illustrated through the Pareto front, and a range of Pareto optimal solutions are shown in Figure \ref{multiexample}(a).

\noindent Figure \ref{multiexample}(b) shows the inventory and energy-related cost values evaluated by the feasible set, aiming like in the GHG case, we visualize and verify the Pareto front. Using lot size ranges $Q_p \in [1000, 1200]$ and $Q_r \in [400, 700]$ with $m = n = s = 1$, we generate $56,079$ combinations of $(C, G)$, displayed as a blue surface. A brute-force algorithm identifies $147$ non-dominated solutions, forming the Pareto front (red curve). Comparing with Figure \ref{multiexample}(a), which uses a coarser grid, confirms that its Pareto points are a subset of those in Figure \ref{multiexample}(b).

\subsubsection*{Existence of Pareto-Optimal Solutions}

\begin{theorem}[Existence of Pareto Solutions]
The bi-objective problems (6) and (9)
admit at least one Pareto-optimal solution.
\end{theorem}

\begin{proof}
\noindent From Theorem \ref{theorem1}, the feasible region is bounded.
Integer variables $(Q_p,Q_r,m,n)$ take finitely many values
within that bounded region.
The objectives $C(\cdot)$, $G(\cdot)$, and $E(\cdot)$
are continuous in $s\in[0,1]$.

Therefore the image of $X$ under the vector objective mapping is a finite union of compact sets. By standard results in multiobjective optimization,
a Pareto-optimal solution exists.
\end{proof}

\medskip

These theoretical results establish that the proposed formulation is a well-posed, bounded, nonconvex mixed-integer nonlinear multiobjective optimization problem with guaranteed existence of efficient solutions.

\subsection{Solver Robustness and Computational Stability}
\noindent Due to the mixed-integer nonlinear structure of the proposed formulation, the optimization problem exhibits nonconvex characteristics that may lead local solvers to converge to suboptimal stationary points depending on initialization and parameter settings. In particular, when solving the scalarized subproblems using BONMIN within AMPL, instances of infeasible convergence and locally optimal solutions were observed for certain weight configurations. To mitigate these effects, multiple initializations and grid refinements were employed to enhance solution stability and improve Pareto front approximation. Furthermore, the brute-force validation approach presented in Figures \ref{multiexampleA} and \ref{multiexample} serves as an independent verification mechanism to ensure that the identified Pareto solutions are not artifacts of solver sensitivity. These findings highlight the computational complexity inherent in multiobjective mixed-integer nonlinear programming problems and suggest that global optimization strategies or hybrid metaheuristic approaches may further enhance robustness in large-scale implementations.

\noindent The purpose of the computational experiments is to demonstrate the applicability and behaviour of the proposed inventory model rather than to compare optimization software packages. BONMIN was selected due to its availability within the AMPL framework and its capability to handle mixed-integer nonlinear programming problems efficiently. Although global optimization solvers such as BARON and SCIP may provide stronger optimality guarantees, the primary objective of this study is to investigate managerial trade-offs arising from the proposed formulation. Comparative evaluation of alternative global optimization approaches remains an interesting direction for future research.

\section{Conclusion}\label{Con}

\noindent In this study, we introduced Economic Production Quantity (EPQ) models within a reverse logistics system. We developed a mathematical formulation to estimate the inventory cost associated with both production and remanufacturing processes. We then formulated a reverse logistics model for EPQ, incorporating constraint \eqref{con1} to ensure sufficient inventory is available throughout the remanufacturing process. In our approach, production and remanufacturing rates are modelled as functions of lot sizes and cycle numbers. Since lot sizes and cycle numbers are integer variables, the problem is formulated as a mixed-integer constrained optimization model. To demonstrate applicability, we tested the model taking examples from existing studies (e.g., Sharma et al. \cite{Swati2021}) and performed sensitivity analysis. Furthermore, to integrate environmental sustainability with economic efficiency, we proposed two extended models that explicitly capture trade-offs between inventory cost and greenhouse gas (GHG) emissions, and between inventory cost and energy consumption, thereby providing an inventory cost function that better reflects real production systems.

\noindent The results, obtained using a variety of solvers in MATLAB and AMPL, demonstrate the model's capability to handle complex inventory systems under sustainability constraints. Numerical experiments (Example \ref{SurfaceEPQ}) show that the proposed model  effectively identifies optimal integer lot sizes for both new and remanufactured items, as well as the corresponding order numbers. In contrast to most models in the literature that treat production and remanufacturing rates as fixed parameters, our formulation models these rates as functions of lot sizes and cycle numbers. This provides a more practical approach, as the rates are closely tied to environmental considerations, such as greenhouse gas emissions and energy consumption. Sensitivity analysis reveals that parameter settings, particularly setup and holding costs, have a significant impact on inventory optimization. Our results provide valuable insights into the effect of such parameters on inventory cost calculation. In addition, in Example \ref{multexam}, we approximated Pareto fronts that contain non-dominated solutions, capturing trade-offs between inventory cost and GHG emissions, and offering decision-makers practical tools to balance economic and environmental objectives. We also approximated Pareto fronts for inventory cost and energy consumption. These results provide managers with valuable guidance for balancing economic efficiency and environmental goals in line with their strategic priorities.

\noindent The present study adopts a deterministic framework to establish the fundamental properties and behaviour of the proposed reverse logistics EPQ model. While deterministic assumptions are commonly employed in the inventory literature to provide analytical insight, uncertainty in demand, return rates, quality recovery, energy consumption, and emission parameters may influence operational decisions in practice. Extending the proposed framework to stochastic, fuzzy, and robust optimization settings represents a promising direction for future research. 

\noindent \textbf{\large {Acknowledgment:}} 
\noindent We acknowledge that during the preparation of this work the author(s) used ChatGPT4.0 in order to rephrase and check grammar. After using this tool/service, the author(s) reviewed and edited the content as needed and take(s) full responsibility for the content of the published article.

\section{\small Funding and/or Conflicts of interests/Competing interests}
 The authors declare that no funds, grants, or other support were received during the preparation of this manuscript. The authors have no relevant financial or non-financial interests to disclose.

 
\section{Appendix}\label{appen}
 
\begin{figure}[H]
	\centering
\includegraphics[scale=1.1]{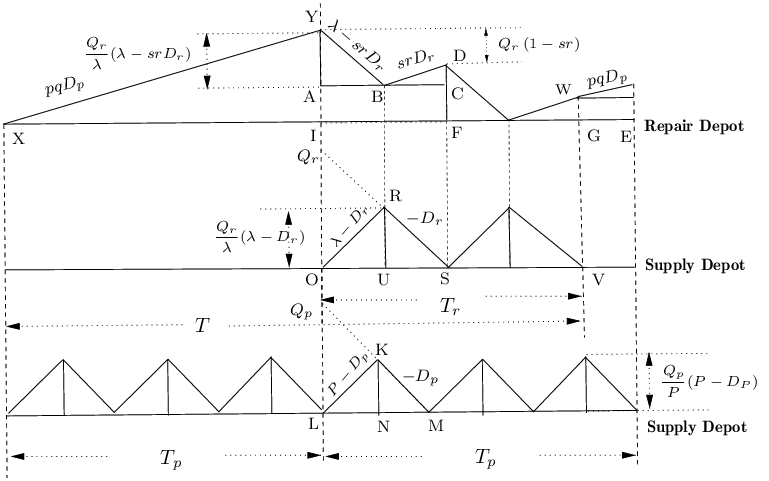}
	\caption{The behaviour of inventory for new, used and repaired items over the interval $T$.}
	\label{figinventory3a}
\end{figure}

\subsection {Inventory cost function introduced in \eqref{eqcost1a}.}
\begin{figure}[H]
	\centering
\includegraphics[scale=1.1]{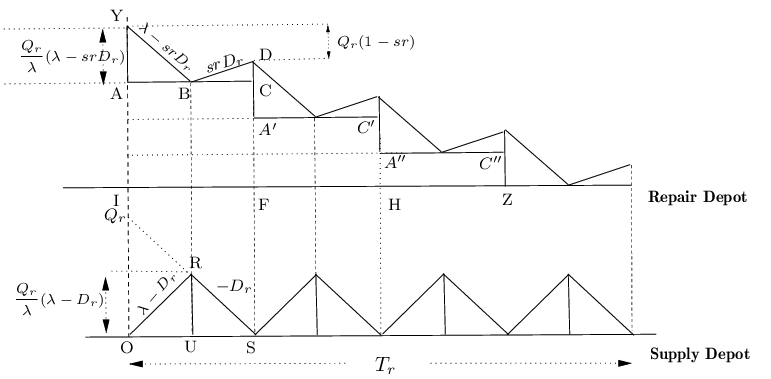}
	\caption{The behaviour of inventory for repaired items over the interval $T_r$.}
	\label{figinventory4}
\end{figure}

\noindent A detailed explanation of the inventory cost function from \eqref{eqcost1a} is now presented, with Figures \ref{figinventory2} and \ref{figinventory4} used to support the corresponding calculations.

The combined area of the supply depot is
\[
n \Delta \mbox{LKM}+m \Delta \mbox{ORS}=\frac{T_pQ_p (P-D_p)}{2P} + \frac{T_rQ_r (\lambda-D_r)}{2 \lambda},
\]
The total area within the repair depot is
\begin{eqnarray*} \label{eqA1}
&=& \Delta \mbox{ XYI}+m\Delta\mbox{YAB}+m\Delta\mbox{BDC}+\rectangle AIFC + \rectangle A^{\prime}FHC^{\prime} + \rectangle A^{\prime\prime}HZC^{\prime\prime} \\
&=& \frac{1}{2}pqD_pT_p^2 + \frac{m}{2}\left(\frac{Q_r}{\lambda}\right)^2(\lambda-srD_r) + \frac{1}{2}srQ_rT_r\left(1-\frac{D_r}{\lambda}\right )^2\\
&+& \left(1-\frac{1}{m}\right)T_r\left(pqD_pT_p-\frac{Q_r}{\lambda}(\lambda - srD_r)-\frac{(m-2)Q_r(1-sr)}{2}\right).
\end{eqnarray*}

Note: Figure \ref{figinventory4} illustrates $\rectangle AIFC$, $\rectangle A^{\prime}FHC^{\prime}$, and $\rectangle A^{\prime\prime}HZC^{\prime\prime}$ for reference.

\noindent \textbf{Author Contributions}: 
Both authors contributed to the study conception and design. Material preparation, data collection and analysis were performed by Dr Indu Bala Wadhawan and Dr Mohammed Mustafa Rizvi.

\end{document}